\documentclass[a4paper,12pt]{article} 
\newdimen\paperhight
\input amssym.def
\input amssym

\newcommand{\ff}{{4\over 5}} 
\newcommand{\st}{{1\over 16}} 
\newcommand{\hf}{\frac{1}{2}}

\newcommand{\qed}{\begin{flushright}{\bf Q.E.D.}\end{flushright}}
\newcommand{\pr}{\par \vspace{3mm} \noindent {\bf [Proof]} \qquad}
\newcommand{\prend}{\hfill \qed}

\newcommand{\e}{\bf e} 
 
\newcommand{\1}{\bf 1}

\newcommand{\be}{\beta}

\newcommand{\la}{\lambda}

\newcommand{\C}{{\Bbb C}} 
\newcommand{\Z}{{\Bbb Z}}

\newcommand{\R}{{\Bbb R}}

\newtheorem{thm}{Theorem}[section]

\newtheorem{lmm}{Lemma}[section]

\begin{document}
\title{3-State Potts model and automorphism of vertex operator algebra 
of order 3}
\author{Masahiko Miyamoto} 
\date{\begin{tabular}{c}
Institute of Mathematics \\
{ University of Tsukuba }\\
{ Tsukuba 305, Japan }\\
\end{tabular}}

\maketitle

\begin{abstract}
In \cite{M2}, the author has defined an automorphism $\tau$ 
of a vertex operator algebra (VOA) of order 2 
using a sub VOA isomorphic to Ising model $L(\hf,0)$.  
We here define an automorphism of VOA of order 3 
by using a sub VOA isomorphic to 
a direct sum of 3-state Potts models $L(\ff,0)$ and an its module 
$L(\ff,3)$. 
If $V$ is the moonshine VOA $V^{\natural}$, 
the defined automorphism is a 3A element of the monster simple group. 
\end{abstract}

\renewcommand{\baselinestretch}{1.2}\large\normalsize
\section{Introduction}
In the research of 
the Griess algebra $V^{\natural}_2$, 
Conway \cite{C} found several idempotents 
called axes of the Griess algebra 
corresponding to elements of the monster simple group 
${\bf M}=Aut(V^{\natural})$.
It was discovered that the Griess algebra 
is the second primary part $(V^{\natural})_2$ 
of the moonshine VOA $V^{\natural}$ constructed in \cite{FLM}. 
  
It was proved in \cite{M2} that idempotents in the Griess algebra 
are halves of conformal vectors (or Virasoro elements of 
sub VOAs). In particular, every idempotent in the Conway's 
list is a half of the Virasoro element of a sub VOA 
isomorphic to one of the minimal 
discrete series of Virasoro VOA $L(n,0)$ with central charge $n$ 
for some $0<n<1$.  
For example, a 2A-involution of 
the monster simple group is a half of the Virasoro element of a sub VOA 
isomorphic to the Ising model $L(\hf,0)$ with central charge $\hf$. 
We note that a sub VOA $W$ in this paper does not usually have the same 
Virasoro element of $V$. 
Conversely, the author showed that if a VOA $V$ contains a sub VOA $W$ 
isomorphic 
to the Ising model $L(\hf,0)$, then it defines an 
automorphism $\tau_W$ of $V$ of order at most 2, which is 
a 2A-element of the monster simple group 
if $V$ is the moonshine VOA $V^{\natural}$. 
The definition is very simple and done as follows: \\
Let $e$ be a Virasoro element of $W$. 
As we will prove the general case in Theorem 5.1, 
$V$ is a direct sum 
of (possibly infinite number of) irreducible $W$-modules since 
$W\cong L(\hf,0)$ is rational. 
It is well known that 
$L(\hf,0)$ has the exactly three irreducible modules $L(\hf,0)$, 
$L(\hf,\hf)$ and $L(\hf,\st)$. We can define an automorphism 
$\tau_e$ by 
$$ \tau_e: \left\{ \begin{array}{rcl}
1& on& \mbox{$W$-submodules isomorphic to } L(\hf,0) \mbox{ or } L(\hf,\hf) \cr
-1& on&\mbox{$W$-submodules isomorphic to } L(\hf,\st) 
\end{array} \right..  $$
Throughout this paper, we will use the similar notation in order 
to define an endomorphism of $V$ by a sub VOA $W$ and  
we will omit "$W$-submodules isomorphic to" from the definition of 
automorphisms in order to simplify 
the notation.  
In the Conway's list, an idempotent for a 3A element is 
a half of the Virasoro element of a sub VOA isomorphic to 
$L(\ff,0)$ with central charge $\ff$. So it is natural for us 
to expect an automorphism $g$ (of order 3) defined by a sub VOA 
isomorphic to $L(\ff,0)$,  
where $L(\ff,0)$ is the third of the discrete series of minimal Virasoro 
vertex operator algebras called 3-state Potts model.  
It is a rational VOA and has the exactly ten irreducible modules : 
$$ \begin{array}{l}
L(\ff,0), L(\ff,\frac{1}{8}),L(\ff,\frac{2}{3}),
L(\ff,\frac{13}{8}),L(\ff,3), \cr
L(\ff,\frac{2}{5}), L(\ff,\frac{1}{40}), 
L(\ff,\frac{1}{15}), L(\ff,\frac{21}{40}), L(\ff,\frac{7}{5}).
\end{array}$$

As we showed in \cite{M2}, if a VOA $V$ contains a sub VOA $W\cong 
L(\ff,0)$, then 
we can define an automorphism $\sigma_W$ of at most 2 given by 
$$ \sigma_W: \left\{ 
\begin{array}{rcl}
1& on& L(\ff,0), L(\ff,3), L(\ff,{13\over 8}), L(\ff,{2\over 5}), 
L(\ff,{1\over 15}), L(\ff,{7\over 5})\cr
-1& on& L(\ff,{1\over 8}), L(\ff,{13\over 8}), L(\ff,{1\over 40)}), 
L(\ff,{21\over 40})  
\end{array} \right. . $$

We note that we can't observe this automorphism in the moonshine VOA. 
Namely, $V^{\natural}$ does not contain any submodules of the second lines. 
Under such a situation, we want to define a 
3A automorphism $\tau_W$ of $V$ 
by a sub VOA $W$ isomorphic to $L(\ff,0)$. 
It is clear that it is not enough to think of only $L(\ff,0)$ 
because there is no difference between the eigenspaces of $\tau$ with 
eigenvalues $e^{2\pi i/3}$ and $e^{4\pi i/3}$. 
Namely, let $V^1$ and $V^2$ 
be eigenspaces of $\tau$ with eigenvalue $e^{2\pi i/3}$ and 
$e^{4\pi i/3}$, 
then they are 
isomorphic as $L(\ff,0)$-modules. 
However, Dong and Mason \cite{DM} showed a wonderful result that 
$V^1$ and $V^2$ are not isomorphic as $V^{<\tau_W>}$-modules, 
where $V^{<\tau_W>}$ is 
the space of $\tau_W$-invariants. So we have to think of a bigger sub VOA.  
What is the difference between 
$V^{<\tau_W>}$ and $L(\ff,0)$?  Recently, Kitazume, Yamada and the author 
have constructed 
a new class of VOAs by using codes over ${\Z}_3$ in \cite{KMY}.  
The interesting point is that 
they used a VOA isomorphic to $L(\ff,0)\oplus L(\ff,3)$ as a sub VOA 
corresponding to $0\in \Z_3$. 

This is our key point 
and the main result in this paper is to show that 
if $V$ contains a sub VOA $W$ isomorphic to 
$L(\ff,0)\oplus L(\ff,3)$, we can define a triality 
automorphism $\tau_W$ of $V$ (of order 3 or possibly 1).

In order to define the automorphism, we need quote 
several results from \cite{KMY}. 
They classified the irreducible modules of $W(0)=L(\ff,0)\oplus L(\ff,3)$.   
Namely,  

\begin{thm}[\cite{KMY}]  
$L(\ff,0)\oplus L(\ff,3)$ is a rational VOA and it has the exactly 
six irreducible modules: 
$$ W(0), W({2\over 5}), 
W({2\over 3},+), W({1\over 15},+),W({2\over 3},-), W({1\over 15},-).$$
Here $h$ in $W(h)$ or $W(h,\pm)$ denotes the lowest degree 
and $W(k,-)$ is a contragredient (dual) module of $W(k,+)$ 
for $k={2\over 3}, {1\over 15}$. 
In particular, 
$$ \begin{array}{l}
W(0)\cong L(\ff,0)\oplus L(\ff,3), \cr
W({2\over 5})\cong L(\ff,{2\over 5})\oplus L(\ff,{7\over 5}),  \cr
W({2\over 3},+)\cong L(\ff,{2\over 3}),  \cr
W({2\over 3},-)\cong L(\ff,{2\over 3}),  \cr
W({1\over 15},+)\cong L(\ff,{1\over 15}),  \cr
W({1\over 15},-)\cong L(\ff,{1\over 15}),  \cr
\end{array} $$
as $L(\ff,0)$-modules. 
\end{thm}

Using the notation in the above theorem, we have:\\

\noindent
{\bf Theorem A}\qquad  
{\it If a VOA $V$ contains a sub VOA $W$ isomorphic to 
$L(\ff,0)\oplus L(\ff,3)$, then 
an endomorphism $\tau_W$ of $V$ defined by 
$$ \tau_W : \left\{ 
\begin{array}{rcl}
1& on& W(0) \mbox{ and } W({2\over 5}) \cr  
e^{2\pi i/3}& on &W({2\over 3},+) \mbox{ and } W({1\over 15},+) \cr
e^{4\pi i/3}& on & W({2\over 3},-) \mbox{ and } W({1\over 15},-)
\end{array} \right.  $$
is an automorphism of $V$. }\\

In order to tell the difference between $W(h,+)$ 
and $W(h,-)$ for $h=\frac{2}{3}, \frac{1}{15}$, 
we have to explain the actions of the lowest 
degree vector of $L(\ff,3)$ since both $W(h,\pm)$ are 
isomorphic to $L(\ff,h)$ as $L(\ff,0)$-modules.  
However, we will take an easier way to avoid such a complicated job.  
We only note that if we fix $W(\frac{2}{3},+)$, 
then $W(\frac{1}{15},\pm)$ and $W(\frac{2}{3},-)$ are 
uniquely determined by the fusion rule 
$$ W({2\over 3},\pm)\times W({2\over 5})=W({1\over 15},\pm). $$ 
We don't distinguish between $W({2\over 3},+)$ and $W({2\over 3},-)$, but 
if we switch them then we shall define $\tau_W^{-1}$.

Let $T$ be a sub VOA of $W$ isomorphic to $L(\ff,0)$. 
If $\tau=1$, then all $T$-submodules of $V$ are 
isomorphic to $L(\ff,0)$, $L(\ff,3)$, $L(\ff,{2\over 5})$ or 
$L(\ff,\frac{7}{5})$.
In this case, we can define another automorphism $\mu_T$ of $V$ as follows: \\

\noindent
{\bf Theorem B} \qquad  
{\it Assume that $V$ contains a sub VOA $T$ isomorphic to $L(\ff,0)$ and 
all $T$-submodules of $V$ are 
isomorphic to $L(\ff,0)$, $L(\ff,3)$, 
$L(\ff,{2\over 5})$ or  $L(\ff,\frac{7}{5})$.
Then the endomorphism $\mu_T$ defined by 
$$ \mu_T : \left\{ 
\begin{array}{rcl}
1& on& L(\ff,0) \mbox{ and } L(\ff,\frac{7}{5}) \cr  
-1& on& L(\ff,3) \mbox{ and } L(\ff,\frac{2}{5}) \cr  
\end{array} \right.  $$
is an automorphism of $V$. } \\

The proofs of these theorems are based on Theorem 2.1 (Proposition 4.4 
in \cite{M1}).  
Namely, 
it is sufficient to show that the fusion rules among the 
irreducible $L(\ff,0)\oplus 
L(\ff,3)$-modules commutes with $\tau_W$.   
For an example, we know the fusion rules among irreducible 
$L(\ff,0)$-modules (Table A), which proves Theorem B. 
Therefore, the main thing we will do in this paper is 
to determine the fusion rules among the irreducible 
$L(\ff,0)\oplus L(\ff,3)$-modules. 

In this paper, we often view $V$ as a $W$-module (or an infinite 
direct sum of $W$-modules) if $V$ contains 
a sub VOA $W$.  
This is not obvious since one of the axioms of VOA-modules expects 
the grade keeping operator $e_1$ of Virasoro element $e$ of $W$ 
to act on $V$ diagonally. 
We will prove in \S 4 that this is generally true 
for a rational sub VOA $W$.

\section{Preliminary results and a generalized VOA constructed 
from a lattice}
Throughout this paper, $W(0)$ denotes a VOA isomorphic to 
$L(\ff,0)\oplus 
L(\ff,3)$. 
Since we will treat only a rational VOA $V$ 
isomorphic to $L(\ff,0)$ or $W(0)$, 
the tensor products of two $V$-modules $M^1$ and $M^2$ 
are always well-defined and it is isomorphic to 
$\oplus_UN^{U}_{M^1,M^2}U, $
where $U$ runs over the all irreducible $V$-modules.
Therefore, it is equal to the fusion rule in our case and so 
we will use the same notation $M^1\times M^2$ to denote the tensor product.

Since $L(\ff,0)\subseteq W(0)$, all $W(0)$-modules are $L(\ff,0)$-modules. 
Using this fact, we will give an upperbound of the fusion rules 
of $W(0)$-modules. 
Using exactly the same proof, we can modify Proposition 11.9 in \cite{DL}
into the following theorem. 

\begin{thm}[\cite{DL}]  
Let $W^1, W^2, W^3$ be $V$-modules  
and assume that $W^1, W^2$ have no proper submodules containing 
$v^1$ and $v^2$, respectively. Let $I\in I\pmatrix{W^3\cr W^1\qquad W^2}$. 
If $I(v^1,z)v^2=0$, then $I(\cdot,z)=0$. 
\end{thm}

In the case where    
$W^1=W(h,+)\oplus W(h,-)$ for $h=\frac{1}{15}, \frac{2}{3}$,  
$W^1$ has no proper submodule containing 
$U^1=\{(v,\phi(v))\in W(h,+)\oplus W(h,-)\}$, where 
$\phi:W(h,+)\to W(h,-)$ is a $L(\ff,0)$-isomorphism.  
Therefore,  we have the following theorem:

\begin{lmm}
The maps 
$$ \begin{array}{l}
\phi^1\ :\ I_{W(0)}\pmatrix{W^3 \cr W(i)
\quad W(j)}
\to I_{L(\ff,0)}\pmatrix{W^3\cr L(\ff,i) \quad L(\ff,j))}, \cr
 \phi^2\ :\ I_{W(0)}\pmatrix{W^3 \cr W(h,+)\oplus W(h,-) 
\quad W(k,+)\oplus W(k,-)}
\to I_{L(\ff,0)}\pmatrix{W^3\cr L(\ff,h) \ \  L(\ff,k)} \cr
\mbox{and}  \cr
 \phi^3\ :\ I_{W(0)}\pmatrix{W^3 \cr W(h,+)\oplus W(h,-) 
\quad W(i)}
\to I_{L(\ff,0)}\pmatrix{W^3\cr L(\ff,h) \quad L(\ff,j)} 
\end{array} $$
induced by the restrictions are all injective 
for $i,j=0,\frac{2}{5}$ and $h,k=\frac{2}{3},\frac{1}{15}$. 
\end{lmm}

Throughout this paper, $N^{W^3}_{W^1, W^2}$ denotes 
$\dim I\pmatrix{W^3\cr W^1\quad W^2}$.

It is known that 
$$
N^{W^3}_{W^1, W^2}=
N^{W^3}_{W^2, W^1}=
N^{(W^1)'}_{W^2, (W^3)'}$$
where $W'$ denotes the contragredient (dual) module of $W$, 
(see Proposition 
5.5.2 in \cite{FHL}). 
We note that $N^{W^1}_{W(0),W^1}=1$ and $N^{W^2}_{W(0),W^1}=0$ 
for $W^1\not\cong W^2$. 
Let $k'=3,{7\over 5}$ for $k=0,{2\over 5}$, respectively. 
Namely, $W(k)\cong L(\ff,k)\oplus L(\ff,k')$ as $L(\ff,0)$-modules. 
By the above lemma and the fusion rules of the 
irreducible $L(\ff,0)$-modules (see Table A),  
we have the following lemma. 

\begin{lmm}
$$\begin{array}{l}
N_{W(i),W(j)}^{W(k)}\leq N_{L(\ff,i),L(\ff,j)}^{L(\ff,k)\oplus L(\ff,k')}\leq 1
\cr
N_{W(i),W(j)}^{W(k,\pm)}\leq N_{L(\ff,i),L(\ff,j)}^{L(\ff,k)}=0 \cr
N_{W(i,+)\oplus W(i,-), W(j,+)\oplus W(j,-)}^{W(k)}
\leq N_{L(\ff,i), L(\ff,j)}^{L(\ff,k)\oplus 
L(\ff,k')}\leq 2 \cr
N_{W(i,+)\oplus W(i,-),W(j,+)\oplus W(j,-)}^{W(k,\pm)}
\leq N_{L(\ff,i),L(\ff,j)}^{L(\ff,k)}\leq 1 
\end{array} $$
\end{lmm}

Let's explain how we determine the fusion rules of $W(0)$-modules. 
Since all VOAs in this paper are rational, we 
identify the fusion product and the tensor product. 
Namely, we will see $W^1\times W^2$ as a $W(0)$-module 
for $W(0)$-modules 
$W^1$ and $W^2$.  Let $T$ be a sub VOA of $W(0)$ isomorphic to $L(\ff,0)$. 
For a $W(0)$-module $W$, it is also a $T$-module and 
we denote them by $(W)_T$. 
As $W(0)$-modules, we have a fusion product $W^1\times W^2$.  
Also viewing $W^1$ and $W^2$ as $L(\ff,0)$-modules, we have 
a fusion product $(W^1)_T\times (W^2)_T$.   
The above two lemmas tell that there is an injective $W(0)$-homomorphism 
$\pi$ of  
$(W^1\times W^2)$ into $(W^1)_T\times (W^2)_T$.  
We note that all fusion rules 
$N^{L(\ff,i)}_{L(\ff,j), L(\ff,h)}$ 
are less than or equal to $1$. 
We will next show 
$N^{W^1}_{W^2, W^3}\not=0$ for desired $W(0)$-modules 
$W^1, W^2, W^3$ so that the above injection $\pi$ of 
$(W^1\times W^2)$ into $(W^1)_T\times (W^2)_T$ is an 
isomorphism.  Therefore, we can 
determine the fusion rules.

In order to show $N^{W^1}_{W^2,W^3}\not=0$ for some $W^1,W^2,W^3$, 
we will use 
a (generalized) VOA $V_L$ constructed from a lattice $L={\sqrt{2}\over 3}A_2$. 
First, we quote a construction of $V_S$ for $S\subseteq {\R}L$ 
for \cite{FLM}. We note that we don't need 
a group extension here since all the values of 
inner products of elements of $L$ are in $2{\Z}[1/3]$.

Let $L$ be a lattice. 
Viewing $H={\C}L$ as a commutative Lie algebra with a bilinear form 
$<,>$, we define the affine Lie algebra 
$$\hat{H}=H[t,t^{-1}]+{\C}k$$ 
associated with $H$ 
and the symmetric tensor algebra $M(1)=S(\hat{H}^-)$ of $\hat{H}^-$, 
where $\hat{H}^-=H[t^{-1}]t$.  As in 
\cite{FLM}, we shall define 
the Fock space $V_L=\oplus_{a\in L}M(1)\e^a$ 
with the vacuum ${\bf 1}=\e^0$ 
and the vertex operators $Y(\ast,z)$ as follows: 

The vertex operator of ${\e}^a$ is given by 
$$ Y(\e^a,z)=exp
\left(\sum_{n\in {\Z}_+}{a(-n)\over n}z^{n}\right)exp
\left(\sum_{n\in {\Z}_+}{a(n)\over -n}z^{-n}\right)
e^{a}z^{a}.  $$
and that of $a(-1)\e^0$ is 
$$ Y(a(-1)\e^0,z)=a(z)=\sum a(n)z^{-n-1}.$$
The vertex operators of other elements are defined by the $n$-th 
normal product:\\
$$  Y(a(n)v,z)=a(z)_nY(v,z)
=Res_x\{(x-z)^na(x)Y(v,z)-(-z+x)^nY(v,z)a(x)\}.  $$
Here the operator of $a\otimes t^n$ on $M(1)\e^b$ are denoted by $a(n)$ and 
$$ \begin{array}{l}
a(n)\e^b=0  \mbox{  for } n>0, \cr
a(0)\e^b=<a,b>\e^b, \cr
e^a\e^b=\e^{a+b}, \cr
z^a\e^b=\e^bz^{<a,b>},\cr
(x+y)^n=\sum_{i=0}^{\infty}{n\choose i}x^{n-i}y^i \mbox{ and }\cr
{n \choose i}={n(n-1)\cdots (n-i+1)\over i!}. 
\end{array}$$
We note that the above definition of vertex operator is very general, that 
is, it is well defined for any $v\in V_{{\R}\otimes L}$ and 
so we may think 
$$  Y(v,z)\in End(V_{\R\otimes L})((z,z^{-1}))   $$
for $v\in V_{\R\otimes_{\Z} L}$, where 
$ P((z,z^{-1}))$ denotes $\{\sum_{n\in \C}a_nz^n: a_n\in P\} $ 
for any set $P$.
In particular, 
$$  Y(v,z)u\in V_{\R\otimes L}[[z,z^{-1}]]z^{\langle a,b\rangle}   $$
for $v\in M(1)e^a$ and $u\in M(1)e^b$ and $a,b\in {\R}L$. 
 Set $\1=\e^0$. 
It is worth to note that if we set $Y(v,z)=\sum_{n\in \R} v_nz^{-n-1}$, 
then $v_{-1}\e^0=v$ for any $v\in \R\otimes L$. 
Also, for the Virasoro element, we set 
$$  w=\sum v^i(-1)^1e^0    $$
where $\{v^1,...,v^k\}$ is an orthonormal basis of ${\R}L$.

For any subset $S$ of ${\R}L$, we can define 
$$ V_S=\oplus_{a\in S}M(1)\e^a. $$ 

The followings are obtained in Chapter 4 of \cite{FLM}.

\begin{lmm}[\cite{FLM}]
$$ \begin{array}{ll}
 [\sum a(n)x^{-n-1}, Y(\e^b,z)]\sim 0 &\mbox{ for any } 
 a,b\in L,  \cr
Y(\e^a,x)Y(\e^b,z)\sim Y(\e^b,z)Y(\e^a,x) &\mbox{ for } <a,b>\equiv 2 \pmod{2},  \cr
Y(\e^a,x)Y(\e^b,z)\sim -Y(\e^b,z)Y(\e^a,z) &\mbox{ for } <a,b>\equiv 1 
\pmod{2}.  
\end{array}  $$
where $a(x,z)\sim b(x,z)$ means  $(z-x)^m(a(x,z)-b(x,z))=0$ for some 
$m\in \Z$. 
Especially, if $\langle a,b\rangle\in 2\Z$, then 
$[Y(v,z),Y(u,x)]\sim 0$ for $v\in M(1)e^a, u\in M(1)e^b$.  
\end{lmm}

In \cite{KMY}, they studied the 
structure of the VOA  $M^0=V_L$ and its modules for 
the lattice $L=\sqrt{2}A_2$.  
Namely, let 
$\langle x,y\rangle=-2$, $\langle x,x\rangle=\langle y,y\rangle=4$ 
and set $L={\Z}x+{\Z}y$ be a lattice (of type $\sqrt{2}A_2$). 
It is easy to see that 
$M^1=V_{{x+2y\over 3}+L}$ and 
$M^2=V_{{2x+y\over 3}+L}$ are $V_L$-modules. 
Set 
$$M=M^0\oplus M^1\oplus M^2.$$  
We note that $M$ is closed 
under the operators $u_n$ of $u\in M$.   
It is proved by \cite{DLMN} that the Virasoro element $w$ of $V_L$ is an 
orthogonal sum of three conformal vectors $w^1$, $w^2$, and $w^3$ with 
central charges $\hf$, ${7\over 10}$, and $\ff$, respectively.
Namely, $V_L$ contains a sub VOA $T$ isomorphic to $L(\hf,0)\otimes 
L(\frac{7}{10},0)\otimes L(\ff,0)$.   
Viewing $V$ as a $T$-module, it is a direct sum of irreducible 
modules of $T$ and each irreducible $T$-module is 
isomorphic to $L(\hf,h_1)\otimes L(\frac{7}{10},h_2)\otimes L(\ff,h_3)$ 
for some $h_1,h_2,h_3$. \\

In the following argument, we recall the study in \cite{KMY}. 

It is clear that $(V_L)_1={\C}x(-1)e^0+{\C}y(-1)e^0$. 
The sum of all subspaces of $M^0=V_{L}$ isomorphic to $L(\hf,0)\otimes 
L(\frac{7}{10},k^1)\otimes L(\ff,k^2)$ for any $k^1, k^2 $ 
is isomorphic to a direct sum of  
$$
T^1=L(\hf,0)\otimes L({7\over 10},0)\otimes \left( L(\ff,0)\oplus L(\ff,3)
\right)$$ and 
$$T^2=L(\hf,0)\otimes L({7\over 10},\frac{3}{5})
\otimes \left( L(\ff,\frac{2}{5})\oplus L(\ff,\frac{7}{5})
\right). $$ 
Also, the sum of subspaces of $M^1$ isomorphic to 
$L(\hf,0)\otimes L(\frac{7}{10},k^1)\otimes L(\ff,k^2)$ 
is a direct sum of  
$$ \begin{array}{l}
T^{3+}=L(\hf,0)\otimes L({7\over 10},0)\otimes L(\ff,{2\over 3}) \cr
T^{4+}=
L(\hf,0)\otimes L({7\over 10},\frac{3}{5})\otimes L(\ff,\frac{1}{15}). 
\end{array} $$
Similarly, 
the sum of subspaces of $M^2$ isomorphic to $L(\hf,0)\otimes 
L(\frac{7}{10},k^1)\otimes L(\ff,k^2)$ is a direct sum of  
$$ \begin{array}{l}
T^{3-}=L(\hf,0)\otimes L({7\over 10},0)\otimes L(\ff,{2\over 3}) \cr
T^{4-}=
L(\hf,0)\otimes L({7\over 10},\frac{3}{5})\otimes L(\ff,\frac{1}{15}). 
\end{array} $$
It is easy to see that 
$T^{n \pm}$ are contragredient (dual) modules of $T^{n \mp}$ 
by the natural inner product of $V_{{\Z}({x+2y\over 3})+L}$.  

Since $\langle L, {x+2y\over 3}+L\rangle\subseteq 2\Z$, 
for $v\in M^{g}$ with $g\in \Z_3$, 
the vertex operator $Y(v,z)$ satisfies $L(-1)$-derivative 
property and the local commutativity 
with all vertex operators $Y(u,z)$ of $u\in V_{L}$. 
By applying it to $M^{h}$ for $h\in \Z_3$, we have 
an intertwining operator 
$Y(\ast,z)\in I\pmatrix{M^{h+g} \cr M^g\qquad M^h}$. 

When we view these intertwining operators as intertwining operators 
among $L(\hf,0)\otimes L(\frac{7}{10},0)\otimes \{
L(\ff,0)\oplus L(\ff,3)\}$-submodules and then 
as intertwining operators among 
$L(\ff,0)\oplus L(\ff,3)$-modules, the following 
theorem is very useful.

\begin{thm}[\cite{M1}] 
Let $W$ be a sub VOA of $V$ which may have a different Virasoro 
element $e\in V_2$.  Assume that $W$ is rational. 
Let $M^1$ and $M^2$ be irreducible $W$-submodules of $V$. 
Set 
$$M(M^1,M^2)=\sum_{v\in M^1, u\in M^2, m\in {\Z}} {\C}v(m)u.$$ 
Then $M(M^1,M^2)$ is a $W$-module and 
we have $I_W{M^3 \choose M^1 \quad M^2}\not= 0$ for any 
irreducible $W$-submodule $M^3$ of $M(M^1,M^2)$.
In particular, let $\{W^i:i\in I\}$ be the set of 
distinct irreducible $W$-modules and 
$V=\oplus V^i$ be the decomposition into 
the direct sum of homogeneous $W$-modules $V^i$, where $V^i$ is 
the sum of all irreducible $W$-submodules of $V$ isomorphic to $W^i$.
For $u\in V$, let $u=\sum u^i$, where $u^i\in V^i$.   
If there are $v\in V^i$, $u\in V^j$, $n\in \Z$ such that 
$(v_nu)^k\not=0$, 
then $I\pmatrix{V^k\cr V^i\qquad V^j}\not=0$.  
\end{thm}

We should note that by Theorem 2.2 and the fusion rules $L(\hf,0)\times 
L(\hf,0)=L(\hf,0)$,
$T^1\oplus T^2\oplus T^{3+}\oplus T^{3-}\oplus T^{4+}\oplus T^{4-}$ 
is closed by the products.

Since $W(h,-)$ and $W(i)$ are contragredient (dual) modules of 
$W(h,+)$ and $W(i)$, respectively, we have:
$$ N^{W(0)}_{W(h,\pm), W(h,\pm)}\not=0 \quad \mbox{ and }\quad 
N^{W(0)}_{W(i), W(i)}\not=0.  \eqno{(2.1)}$$
for $h=0,\frac{2}{5}$ and $i=\frac{2}{3},\frac{1}{15}$. 

It is easy to check that 
$$\begin{array}{l}
(T^2)_1={\C}(x+2y)(-1)e^0 \cr
(T^1)_2={\C}(3x(-1)^2e^0+(x+2y)(-1)^2e^0) \cr 
(T^{3+})_{2/3}={\C}(e^{(x+2y)/3}+e^{(x-y)/3}+e^{(-2x-y)/3}) \cr
(T^{4+})_{2/3}={\C}(2e^{(x+2y)/3}-e^{(x-y)/3}-e^{(-2x-y)/3}) \cr
\end{array} $$
Since  
$$ ((x+2y)(-1)e^0)_{-1}(x+2y)(-1)e^0=(x+2y)(-1)^2e^0\not\in (T^1)_2 $$
we have 
$$ N^{W(\frac{2}{5})}_{W(\frac{2}{5}), W(\frac{2}{5})}\not=0. 
\eqno{(2.2)}$$
Also, since 
$$(x+2y)(-1)\{\lambda e^{(x+2y)/3}+\mu (e^{(x-y)/3}+e^{(-2x-y)/3})\}
=4\lambda e^{(x+2y)/3}-2\mu(e^{(x+2y)/3}+e^{(-2x-y)/3}), $$
we have :
$$
N^{W(\frac{1}{15},+)}_{W(\frac{2}{5}), W(\frac{2}{3},+)}= 
N^{W(\frac{2}{5})}_{W(\frac{2}{3},+), W(\frac{1}{15},-)}
=N^{W(\frac{2}{3},-)}_{W(\frac{2}{5}), W(\frac{1}{15},-)}\not=0 \eqno{(2.3)} $$
and 
$$
N^{W(\frac{1}{15},+)}_{W(\frac{2}{5}), 
W(\frac{1}{15},+)}=
N^{W(\frac{2}{5})}_{W(\frac{1}{15},+), 
W(\frac{1}{15},-)}\not=0. \eqno{(2.4)} $$
Similarly, 
$$
N^{W(\frac{1}{15},-)}_{W(\frac{2}{5}), 
W(\frac{2}{3},-)}= 
N^{W(\frac{2}{5})}_{W(\frac{2}{3},-), 
W(\frac{1}{15},+)}
=N^{W(\frac{2}{3},+)}_{W(\frac{2}{5}), 
W(\frac{1}{15},+)}\not=0 \eqno{(2.5)} $$
and 
$$ N^{W(\frac{1}{15},-)}_{W(\frac{2}{5}), 
W(\frac{1}{15},-)}= 
N^{W(\frac{2}{5})}_{W(\frac{1}{15},-), 
W(\frac{1}{15},+)}\not=0. \eqno{(2.6)} $$

It follows from the direct calculations that 
$e^{(x+2y)/3}+e^{(x-y)/3}+e^{(-2x-y)/3}$ 
and 
$2e^{(x+2y)/3}-(e^{(x-y)/3}+e^{(-2x-y)/3})$ 
are  
lowest degree vectors of 
$L(\hf,0)\otimes 
L(\frac{7}{10},0)\otimes L(\ff,\frac{2}{3})\subseteq M^1$ and   
$L(\hf,0)\otimes 
L(\frac{7}{10},\frac{3}{5})\otimes L(\ff,\frac{1}{15})\subseteq M^1$, 
respectively. Similarly,
$e^{(-x-2y)/3}+e^{(-x+y)/3}+e^{(2x+y)/3}$ 
and 
$2e^{(-x-2y)/3}-(e^{(-x+y)/3}+e^{(2x+y)/3})$ 
are  
lowest degree vectors of 
$L(\hf,0)\otimes 
L(\frac{7}{10},0)\otimes L(\ff,\frac{2}{3})\subseteq M^2$ and   
$L(\hf,0)\otimes 
L(\frac{7}{10},\frac{3}{5})\otimes L(\ff,\frac{1}{15})\subseteq M^2$, 
respectively. 
Also, for 
$u={\alpha}e^{(x+2y)/3}+{\be}(e^{(x-y)/3}+e^{(-2x-y)/3})$ and 
$v={\lambda}e^{(x+2y)/3}+{\mu}(e^{(x-y)/3}+e^{(-2x-y)/3})$, 
we have $u_{-1/3}v={2\be\mu}e^{(x+2y)/3}+{\be\lambda+\alpha\mu}
(e^{(x-y)/3}+e^{(-2x-y)/3})$, 
where $u_{-1/3}$ is the grade keeping operator of $u$.
Hence, we have  
$$
N^{W(\frac{2}{3},\mp)}_{W(\frac{2}{3},\pm), 
W(\frac{2}{3},\pm)}\not=0, \eqno{(2.7)} $$
$$ N^{W(\frac{2}{3},\mp)}_{W(\frac{1}{15},\pm), 
W(\frac{1}{15},\pm)}=
N^{W(\frac{1}{15},\mp)}_{W(\frac{1}{15},\pm), W(\frac{2}{3},\pm)}
\not=0, \eqno{(2.8)} $$
$$N^{W(\frac{1}{15},\mp)}_{W(\frac{1}{15},\pm), W(\frac{1}{15},\pm)}
\not=0. \eqno{(2.9)}$$

\section{Fusion rule}
We first list the fusion rules among $L(\ff,0)$-modules 
$L(\ff,0)$, $L(\ff,3)$, $L(\ff,2/5)$, $L(\ff,7/5)$, 
$L(\ff,\frac{2}{3})$ and $L(\ff,\frac{1}{15})$, which are the only 
irreducible modules we need in this section. 
For the fusion rules for the remaining cases, see Appendix. 

\begin{center}
{Table A}
\end{center}
$$ \begin{array}{|c||c|c|c|c|c|}
\hline
0& \frac{2}{5}& \frac{7}{5}& \frac{1}{15}& 3& \frac{2}{3}  \cr
\hline 
\hline
\frac{2}{5}& 0\!:\!\frac{7}{5}& \frac{2}{5}\!:\!3& \frac{1}{15}\!:\!\frac{2}{3}& \frac{7}{5}& \frac{1}{15}  \cr
\hline 
\frac{7}{5}& \frac{2}{5}\!:\!3& 0\!:\!\frac{7}{5}& \frac{2}{3}\!:\!\frac{1}{15}& \frac{2}{5}
& \frac{1}{15} \cr
\hline 
\frac{1}{15}& \frac{1}{15}\!:\!\frac{2}{3}& \frac{2}{3}\!:\!\frac{1}{15}& 
0\!:\!\frac{7}{5}\!:\!\frac{2}{3}\!:\!\frac{1}{15}\!:\!3\!:\!\frac{2}{5}& \frac{1}{15}& \frac{2}{5}\!:\!\frac{1}{15}\!:\!\frac{7}{5}  \cr
\hline 
3& \frac{7}{5}& \frac{2}{5}& \frac{1}{15}& 0& \frac{2}{3}  \cr
\hline 
\frac{2}{3}& \frac{1}{15}& \frac{1}{15}& \frac{2}{5}\!:\!\frac{1}{15}\!:\!\frac{7}{5}& \frac{2}{3} 
& 0\!:\!\frac{2}{3}\!:\!3  \cr
\hline 
\end{array} $$
In the table, the number $h$ denotes $L(\ff,h)$ and 
$h:\cdots :k$ denotes $L(\ff,h)+\cdots +L(\ff,k)$.

By (2.5) and (2.6), 
$N^{W(h,\mp)}_{W(h,\pm), W(h,\pm)}\not=0$. 
Hence, by the fusion rules 
$L(\ff,{2\over 3})\times L(\ff,{2\over 3})
=L(\ff,0)+L(\ff,3)+L(\ff,{2\over 3})$
of $L(\ff,0)$-modules and (2.1) and (2.7), 
we have 
$$ \begin{array}{l}
W({2\over 3},\pm)\times W({2\over 3},\pm)=W({2\over 3},\mp)\cr
W({2\over 3},\pm)\times W({2\over 3},\mp)=W(0). \end{array}
\eqno{(3.1)}  $$

Similarly, by the fusion rules 
$L(\ff,{1\over 15})\times L(\ff,{1\over 15})
=L(\ff,0)+L(\ff,3)+L(\ff,{2\over 3})+L(\ff,{1\over 15})+L(\ff,{2\over 5})
+L(\ff,{7\over 5})$
of $L(\ff,0)$-modules and (2.1),(2.4),(2.8) and (2.9), 
we have 
$$\begin{array}{l}
W({1\over 15},\pm)\times W({1\over 15},\pm)
=W({1\over 15},\mp)+W(\frac{2}{3},\mp) \cr
W({1\over 15},\pm)\times W({1\over 15},\mp)
=W(0)+W(\frac{2}{5})\end{array}
 \eqno{(3.2)} $$
 
By the fusion rules 
$L(\ff,{2\over 5})\times L(\ff,{2\over 3})
=L(\ff,{1\over 15})$ and (2.3) and (2.5), 
we have 
$$ W(\frac{2}{5})\times W(\frac{2}{3},\pm)=W(\frac{1}{15},\pm). \eqno{(3.3)}$$
By the fusion rules 
$L(\ff,{2\over 5})\times L(\ff,{1\over 15})
=L(\ff,{1\over 15})+L(\ff,\frac{2}{3})$ and (2.3) $\sim$ (2.6), 
we have 
$$ W(\frac{2}{5})\times W(\frac{1}{15},\pm)=W(\frac{1}{15},\pm)+
W(\frac{2}{3},\pm). \eqno{(3.4)}$$
By the fusion rules 
$L(\ff,{2\over 5})\times L(\ff,{2\over 5})
=L(\ff,0)+L(\ff,\frac{7}{5})$ and (2.1) and (2.2), 
we have 
$$ W(\frac{2}{5})\times W(\frac{2}{5})=W(0)+W(\frac{2}{5}). \eqno{(3.5)}$$
By the fusion rules 
$L(\ff,{2\over 3})\times L(\ff,{1\over 15})
=L(\ff,\frac{2}{5})+L(\ff,\frac{7}{5})+L(\ff,\frac{1}{15})$ 
and (2.3),(2.5) and (2.8), 
we have 
$$ \begin{array}{l}
W(\frac{2}{3},\pm)\times W(\frac{1}{15},\pm)
=W(\frac{1}{15},\mp) \cr
W(\frac{2}{3},\pm)\times W(\frac{1}{15},\mp)
=W(\frac{2}{5}).  \end{array}
 \eqno{(3.6)}$$

We put the above fusion rules in the following table. 

\begin{center}
{Table B}
\end{center}
\small
$$ \begin{array}{|c||c|c|c|c|c|}
\hline
W(0) \! &\!  W(\frac{2}{5}) \! &\!  W(\frac{2}{3},+) \! &\! W(\frac{1}{15},+) \! &\! 
W(\frac{2}{3},-) \! &\! W(\frac{1}{15}.-) \cr
\hline 
\hline
W(\frac{2}{5})\! &\! W(0)\!:\! W(\frac{2}{5})\! &\! W(\frac{1}{15},+)\! &\!  W(\frac{1}{15},+)\!:\! 
W(\frac{2}{3},+) \! &\!  W(\frac{1}{15},-)\! &\!  W(\frac{1}{15},-)\!:\! 
W(\frac{2}{3},-) \cr
\hline
W(\frac{2}{3},+) \! &\! W(\frac{1}{15},+)\! &\! W(\frac{2}{3},-)\! &\! W(\frac{1}{15},-)\! &\! W(0)\! &\! W(\frac{2}{5})   \cr
\hline
W(\frac{1}{15},+) \! &\!   W(\frac{1}{15},+)\!:\! 
W(\frac{2}{3},+)\! &\!  W(\frac{1}{15},-)\! &\!  W(\frac{1}{15},-)\!:\! 
W(\frac{2}{3},-) \! &\!  W(\frac{2}{5}) \! &\!  W(0)\!:\! W(\frac{2}{5})\cr
\hline
W(\frac{2}{3},-)\! &\! W(\frac{1}{15},-)\! &\!  W(0)\! &\!  W(\frac{2}{5}) \! &\! W(\frac{2}{3},+) \! &\! W(\frac{1}{15},+) \cr
\hline
W(\frac{1}{15},-)\! &\!   W(\frac{1}{15},-)\!:\! 
W(\frac{2}{3},-)\! &\! W(\frac{2}{5})\! &\!  W(0)\!:\! W(\frac{2}{5}) \! &\! W(\frac{1}{15},+)\! &\!  W(\frac{1}{15},+)\!:\! 
W(\frac{2}{3},+) \cr
\hline
\end{array} $$
\normalsize

\section{Automorphisms}
As we showed in \cite{M2}, if a VOA contains $L(\ff,0)$, then 
we have an automorphism $\sigma$ of at most 2 given by 
$$ \sigma: \left\{ 
\begin{array}{rcl}
1& on& L(\ff,0), L(\ff,3), L(\ff,{2\over 3}), L(\ff,{2\over 5}), 
L(\ff,{1\over 15}), L(\ff,{7\over 5})\cr
-1& on& L(\ff,{1\over 8}),L(\ff,{13\over 8}), L(\ff,{1\over 40}), 
L(\ff,{21\over 40}) 
\end{array} \right. . $$

So we next think about the case $\sigma=1$ or the space $V^{\sigma}$ of 
$\sigma$-invariants. 
In this case, there are no $L(\ff,0)$-submodules isomorphic to 
$L(\ff,{1\over 8}), L(\ff,{13\over 8}), L(\ff,{1\over 40})$ or  
$L(\ff,{21\over 40})$. We next assume that $V$ contains 
$L(\ff,0)\oplus L(\ff,3)$. We should note that if $V$ 
contains $L(\ff,0)\oplus L(\ff,3)$, then 
there are no $L(\ff,0)$-submodules in $V$ isomorphic to 
$L(\ff,{1\over 8}), L(\ff,{13\over 8}), L(\ff,{1\over 40})$ or  
$L(\ff,{21\over 40})$. \\

\noindent
{\bf Theorem A} \qquad 
{\it If a VOA $V$ contains $L(\ff,0)\oplus L(\ff,3)$, then 
an endomorphism $\tau$ of $V$ defined by 
$$ \tau : \left\{ 
\begin{array}{rcl}
1\ & on& W(0) \mbox{ and } W({2\over 5}) \cr  
e^{2\pi i/3}\ & on &W({2\over 3},+) \mbox{ and } W({1\over 15},+) \cr
e^{4\pi i/3}\ & on & W({2\over 3},-) \mbox{ and } W({1\over 15},-)
\end{array} \right.  $$
is an automorphism of $V$. }\\

\pr 
Replacing $W(i)$ and $W(h,+)$ and $W(k,-)$ in the table (B) by 
$1$ and $e^{2\pi i/3}$ and $e^{4\pi i/3}$, we have
$$ \begin{array}{|c||c|c|c|c|c|}
\hline
1 & 1 & e^{2\pi i/3} &e^{2\pi i/3} &
e^{4\pi i/3} &e^{4\pi i/3} \cr
\hline 
\hline
1&1:1&e^{2\pi i/3}& e^{2\pi i/3}:
e^{2\pi i/3} & e^{4\pi i/3}& e^{4\pi i/3}:
e^{4\pi i/3} \cr
\hline
e^{2\pi i/3} &e^{2\pi i/3}&e^{4\pi i/3}&e^{4\pi i/3}&1&1   \cr
\hline
e^{2\pi i/3} &  e^{2\pi i/3}:
e^{2\pi i/3}& e^{4\pi i/3}& e^{4\pi i/3}:
e^{4\pi i/3} & 1 & 1:1\cr
\hline
e^{4\pi i/3}&e^{4\pi i/3}& 1& 1 &e^{2\pi i/3} &e^{2\pi i/3} \cr
\hline
e^{4\pi i/3}&  e^{4\pi i/3}:
e^{4\pi i/3}&1& 1:1 &e^{2\pi i/3}& e^{2\pi i/3}:
e^{2\pi i/3} \cr
\hline
\end{array} $$
which is compatible with the products. Hence, by Theorem 2.2, 
$\tau$ is an automorphism of $V$. 
\prend

If $\tau_W=1$, then all $T$-submodules of $V$ are 
isomorphic to $L(\ff,0)$, $L(\ff,3)$, $L(\ff,{2\over 5})$ or 
$L(\ff,\frac{7}{5})$ for $T\subseteq W$ and $T\cong L(\ff,0)$.
In this case, we can define another automorphism $\mu_T$ of $V$ as follows: \\

\noindent
{\bf Theorem B} \qquad  
{\it Assume that $V$ contains a sub VOA $T$ isomorphic to $L(\ff,0)$ and 
all $T$-submodules of $V$ are 
isomorphic to $L(\ff,0)$, $L(\ff,3)$, 
$L(\ff,{2\over 5})$ or  $L(\ff,\frac{7}{5})$.
Then the endomorphism $\mu_T$ defined by 
$$ \mu_T : \left\{ 
\begin{array}{rcl}
1& on& L(\ff,0) \mbox{ and } L(\ff,\frac{7}{5}) \cr  
-1& on& L(\ff,3) \mbox{ and } L(\ff,\frac{2}{5}) \cr  
\end{array} \right.  $$
is an automorphism of $V$. } \\

\section{$V$ as a sub VOA-module}
The notion of sub VOAs of $V$ in this paper is not the same as 
in \cite{FZ}, where 
they expected sub VOA $W$ to have the same Virasoro element with $V$. 
Our definition of sub VOAs is:  \\
$(W=\oplus W_n, Y^W, w_W,\1_W)$ is a sub VOA of $V$ if \\
(1) $(W,Y^W,w_W,\1_W)$ is a VOA, \\
(2) $W\subseteq V$ and $W_n=W\cap V_n$, \\
(3) $\1_W=\1_V$ and  \\
(4) $Y^W(v,z)=Y(v,z)|_W$ for $v\in W$. \\

There are several definitions for VOA-modules, but   
we will include an infinite direct sum of irreducible modules as a 
VOA-module $M$.  Namely, we don't assume $\dim M_n<\infty$.

Let $W$ be a sub VOA of $V$ and $e$ a Virasoro element of $W$.  
Different from the ordinary algebras, it is not obvious that $V$ 
is a $W$-module. The problem is whether $e_1$ acts on $V$ 
diagonally or not. 

The purpose of this section is to show that 
$V$ is a $W$-module for $W$ in our cases. 

Let $V$ be a VOA and $W$ a sub VOA. 
Let $w$ and $e$ be Virasoro elements of $V$ and $W$, respectively. 

In particular, $e\in V_2$ and so $e_1$ keeps the grade of $V$. 
By the assumption, $f_1=w_1-e_1$ acts on $W$ as $0$. 
Furthermore, for $v\in W$ 
$Y(e_0v,z)|_W={d\over dz}Y(v,z)|_W=Y(w_0v,z)|_W$ and so 
$Y(w_0v-e_0v,z)|_W=0$. 
In particular, $w_0v-e_0v=(w_0v-e_0v)_{-1}{\bf 1}=0$ and so, 
$f_0=w_0-e_0$ acts on $W$ as $0$. Hence we have:

\begin{lmm}  Both $f_1$ and $f_0$ commutes $v_m$ on $V$ 
for $v\in W$ and $m\in \Z$.
\end{lmm}

\pr
It follows from $[f_0,v_m]=(f_0v)_m=0 $ and 
$[f_1,v_m]=(f_0v)_{m+1}+(f_1v)_m=0 $.
\prend

The main purpose in this section is to prove the following theorem:

\begin{thm} If $W$ is rational, then $V$ is a $W$-module. 
\end{thm}

\pr 
Define a module vertex operator $Y^V(v,z)$ of $v\in W$ by 
the vertex operator of $v\in V$. 
Clearly, they satisfy the local commutativity and 
the $e_0$-derivative property:
$$Y^V(e_0v,z)=Y^V((f_0+e_0)v,z)=Y(w_0v,z)={d\over dz}Y(v,z) $$
for $v\in W$. 
Hence, what we have to do is to 
prove that $V$ is a direct sum of eigenspaces of $e_1$. 
Suppose false. 
Since $f_1$ commutes with all $v_n$ for $v\in W$, 
the eigenspace $V_{\la}$  and the generalized eigenspace 
$T_{\la}=\{v\in V| \exists n\in {\Z}\quad (f_1-\la)^nv=0 \}$ of 
$f_1$ with eigenvalue $\la$ is invariant 
under the actions $v_n$ of $v\in W$. 
We first prove that the eigenspace $V_{\la}$ is 
a direct sum of irreducible $W$-modules. 
Since $f_1$ keeps the grade, it acts on $V_n$ and so we have 
$V_{\la}=\oplus (V_{\la})_n$, where 
$(V_{\la})_n=V_n \cap V_{\la}$. Since $e_1=(w_1-f_1)$ 
acts on $(V_{\la})_n$ as $n-\la$, $V_{\la}$ is a $W$-module. 
Since $W$ is rational by the assumption, $V_{\la}$ is a 
direct sum of irreducible $W$-modules. 

Since $f_1$ acts on each finite dimensional homogenous spaces $V_n$, 
$V$ is a direct sum of generalized eigenspaces of $f_1$. 
Hence, there are $\la$ and $n$ such that 
$T_{\la}\cap V_n\not=V_{\la}\cap V_n$.  Take $n$ as a minimal one. 

As we explained as above, $v_n$ acts on 
$T_{\la}/V_{\la}$ for $v\in W$ and 
the eigenspace $X_{\la}$ of $f_1$ in 
$T_{\la}/V_{\la}$ is a $W$-module.  By the choice of 
$n$,  
$(X_{\la})_n\not=0$. 
Let $\bar{X}$ be an irreducible $W$-submodule of $X_{\la}$ 
whose lowest degree is $n$.
We should note that since $V_{\la}$ is an eigenspace of $f_1$ 
and $e_1=w_1-f_1$, 
the lowest eigenspace of $w_1$ in $V_{\la}$ is the lowest 
eigenspace  of $e_1$.

Hence there is an irreducible $W$-submodule $\bar{B}$ of $T_{\la}/V_{\la}$ 
whose lowest degree space $\bar{B}_0$ is 
in $(V_n+V_{\la})/V_{\la}$ since $W$ is 
rational.
Let $B$ be its inverse image. Clearly, $B$ contains $V_{\la}$ and 
$f_1$ does not act on $B$ diagonally. 
Let $S$ is a submodule of $V_{\la}$ generated 
by $W$-submodules which are not isomorphic to $\bar{B}$. 
Since $(f_1-\la)B\not=0$ and all submodule of 
$(f_1-\la)B$ is isomorphic to $\bar{B}$, all composition factors of 
$B/S$ is isomorphic to $\bar{B}$ as $W$-modules and 
we have $S\not= V_{\la}$. In particular, 
$(B/S)\cap (V_m+S/S)=0$ for all $m<n$. 
We next show that Zhu-algebra $A(W)$ acts on the top module 
$(B_n+S)/S$ 
of $B/S$.   

In order to prove the above assertion, 
we will use an idea for Zhu-algebra in \cite{Z}.
We will treat a general case for a while. 
Let $A(V)=V/O(V)$ be the Zhu-algebra of $V$.
For $v\in V$, $o(v)$ denotes the grade keeping operator of $v$.
For a homogeneous element $v\in V_m$, 
if $Y^M(v,z)=\sum v_iz^{-i-1}$ is a module vertex operator, then 
$o(v)=v_{m-1}$.  It actually depends on the module $M$, 
but we write $o(v)=v_{m-1}$ formally. 

Let $R$ be the ring generated by all $o(v)$ of $v\in V$. 
Let $I$ be a subspace of $R$ generated by 
the elements $\sigma=v^1_{i_1}\cdots v^r_{i_r}\in R$ satisfying 
that $v^s_{i_s}\cdots v^r_{i_r}$ decreases the grade for some 
$1\leq s\leq r$. 
We permit an infinite sum of such elements if it is well-defined in $R$.
Clearly, $I$ is a two-sided ideal of $R$.  
It is known that $A(V)=R/I$, see \cite{Z}.

Let's go back to the proof. 
Since $W$ is rational, 
$A(W)$ is a semi-simple.  
Let $\phi=v^1_{i_1}\cdots v^t_{i_t}\in R$ and assume that 
$v^s_{i_s}\cdots v^t_{i_t}$ decreases the grade on $W$. 
Since the grade on $W$ is the same as that on $V$,    
$\phi$ acts on 
$(B_n+S)/S$ trivially.
Hence $A(W)$ acts on $(B_n+S)/S$.  Since $A(W)$ is semi-simple, 
$(B_n+S)/S$ is a direct sum of irreducible $A(W)$-modules. 
By the definition of $B$ and $S$, 
$(B_n+S)/S$ is a homogeneous $A(W)$-module. 
Since $e_1$ is in the center of $A(W)$, 
$e_1$ acts on $(B_n+S)/S$ as a scalar times and so does $f_1$,
which contradicts to the facts that 
$(f_1-\la)(B_n+S)/S\not=0$. 

Hence, $V$ is a direct sum of eigenspaces of $f_1$ and so 
$e_1$ acts on $V$ diagonally. 
This completes the proof of Theorem 5.1. 
\prend

\section{Appendix} 
\subsection{Fusion rules of irreducible $L(\ff,0)$-modules}
For $L(\ff,0)$, the following fusion rules are known, see \cite{W}. 
In the following table, the numbers $h$ denote $L(\ff,h)$ and 
$h_1:...:h_t$ denotes $L(\ff,h_1)+...+h(\ff,h_t)$.

\begin{center}
{Table C}
\end{center}
\small
$$ \begin{array}{|c||c|c|c|c|c|c|c|c|c|}
\hline
0& \frac{2}{5}& \frac{1}{40}& \frac{7}{5}& \frac{21}{40}& \frac{1}{15}& 3& 
\frac{13}{8}& \frac{2}{3}& \frac{1}{8}  \cr
\hline 
\hline
\frac{2}{5}& 0\!:\!\frac{7}{5}& \frac{1}{8}\!:\!\frac{21}{40}& \frac{2}{5}\!:\!3& \frac{1}{40}\!:\!\frac{13}{8}& \frac{1}{15}\!:\!\frac{2}{3}& \frac{7}{5}& \frac{21}{40}& \frac{1}{15}& \frac{1}{40}  \cr
\hline 
\frac{1}{40}& \frac{1}{8}\!:\!\frac{21}{40}& 0\!:\!\frac{7}{5}\!:\!\frac{2}{3}\!:\!\frac{1}{15}& \frac{1}{40}\!:\!\frac{13}{8}& \frac{2}{5}\!:\!3\!:\!\frac{1}{15}\!:\!\frac{2}{3}& 
\frac{1}{40}\!:\!\frac{13}{8}\!:\!\frac{21}{40}\!:\!\frac{1}{8}& \frac{21}{40}& \frac{7}{5}\!:\!\frac{1}{15}& \frac{21}{40}\!:\!\frac{1}{40}& \frac{1}{15}\!:\!\frac{2}{5}  \cr
\hline 
\frac{7}{5}& \frac{2}{5}\!:\!3& \frac{1}{40}\!:\!\frac{13}{8}& 0\!:\!\frac{7}{5}& \frac{1}{8}\!:\!\frac{21}{40}& \frac{2}{3}\!:\!\frac{1}{15}& \frac{2}{5}
& \frac{1}{40}& \frac{1}{15}& \frac{21}{40} \cr
\hline 
\frac{21}{40}& \frac{1}{40}\!:\!\frac{13}{8}& \frac{2}{5}\!:\!3\!:\!\frac{1}{15}\!:\!\frac{2}{3}& \frac{1}{8}\!:\!\frac{21}{40}& 0\!:\!\frac{7}{5}\!:\!\frac{2}{3}\!:\!\frac{1}{15}& \frac{1}{8}\!:\!\frac{21}{40}
\!:\!\frac{13}{8}\!:\!\frac{1}{40}& \frac{1}{40}& \frac{2}{5}\!:\!\frac{1}{15}& \frac{1}{40}\!:\!\frac{21}{40}& \frac{1}{15}\!:\!\frac{7}{5}  \cr
\hline 
\frac{1}{15}& \frac{1}{15}\!:\!\frac{2}{3}& \frac{1}{40}\!:\!\frac{13}{8}\!:\!\frac{21}{40}\!:\!\frac{1}{8}& \frac{2}{3}\!:\!\frac{1}{15}& \frac{1}{8}\!:\!\frac{21}{40}\!:\!\frac{13}{8}\!:\!\frac{1}{40}& 
0\!:\!\frac{7}{5}\!:\!\frac{2}{3}\!:\!\frac{1}{15}\!:\!3\!:\!\frac{2}{5}& \frac{1}{15}& \frac{1}{40}\!:\!\frac{21}{40}& \frac{2}{5}\!:\!\frac{1}{15}\!:\!\frac{7}{5}& \frac{1}{40}\!:\!\frac{21}{40}  \cr
\hline 
3& \frac{7}{5}& \frac{21}{40}& \frac{2}{5}& \frac{1}{40}& \frac{1}{15}& 0& \frac{1}{8}& \frac{2}{3}& \frac{13}{8}  \cr
\hline 
\frac{13}{8}& \frac{21}{40}& \frac{7}{5}\!:\!\frac{1}{15}& \frac{1}{40}& \frac{2}{5}\!:\!\frac{1}{15}& \frac{1}{40}\!:\!\frac{21}{40}& \frac{1}{8}& 0\!:\!\frac{2}{3} 
& \frac{1}{8}\!:\!\frac{13}{8}& \frac{2}{3}\!:\!3  \cr
\hline 
\frac{2}{3}& \frac{1}{15}& \frac{21}{40}\!:\!\frac{1}{40}& \frac{1}{15}& \frac{1}{40}\!:\!\frac{21}{40}& \frac{2}{5}\!:\!\frac{1}{15}\!:\!\frac{7}{5}& \frac{2}{3} 
& \frac{1}{8}\!:\!\frac{13}{8}& 0\!:\!\frac{2}{3}\!:\!3& \frac{1}{8}\!:\!\frac{13}{8}  \cr
\hline 
\frac{1}{8}& \frac{1}{40}& \frac{1}{15}\!:\!\frac{2}{5}& \frac{21}{40}& \frac{1}{15}\!:\!\frac{7}{5}& \frac{1}{40}\!:\!\frac{21}{40}& \frac{13}{8}& \frac{2}{3}\!:\!3
& \frac{1}{8}\!:\!\frac{13}{8}& 0\!:\!\frac{2}{3} \cr
\hline 
\end{array} $$

\end{document}